\documentclass[11pt,a4paper]{article}
\usepackage[T1]{fontenc}
\usepackage[utf8]{inputenc}
\usepackage{lmodern}
\usepackage{amsmath,amssymb,amsthm,mathtools}
\usepackage[margin=1.13in,top=1.03in,bottom=1.06in,headheight=14pt]{geometry}
\usepackage[protrusion=true,expansion=true]{microtype}
\usepackage{enumitem}
\usepackage{xcolor}
\usepackage{fancyhdr}
\usepackage[bookmarksnumbered,unicode]{hyperref}
\usepackage{bookmark}
\usepackage{etoolbox}
\definecolor{inkblue}{HTML}{204E64}
\hypersetup{colorlinks=true,linkcolor=inkblue,citecolor=inkblue,urlcolor=inkblue,
 pdfauthor={Ma Yicen},pdftitle={AP without BAP in a separable Banach space and its dual},
 pdfsubject={AP-BAP conclusion obtained October 30, 2026; revised November 5, 2026},
 pdfkeywords={approximation property, bounded approximation property, dual Banach space}}
\newcommand{\doi}[1]{\mbox{\href{https://doi.org/#1}{doi:\nolinkurl{#1}}}}
\newcommand{\arxiv}[1]{\mbox{\href{https://arxiv.org/abs/#1}{arXiv:\nolinkurl{#1}}}}
\setlist[enumerate]{leftmargin=1.8em,itemsep=2pt,topsep=4pt}
\newtheorem{theorem}{Theorem}[section]
\newtheorem{lemma}[theorem]{Lemma}
\newtheorem{proposition}[theorem]{Proposition}
\newtheorem{corollary}[theorem]{Corollary}
\theoremstyle{definition}

\theoremstyle{remark}

\begin{document}
\begin{center}
{\LARGE\bfseries AP without BAP in a separable Banach space}\\[3pt]
{\LARGE\bfseries and its dual}\\[8pt]
{\normalsize\itshape A direct-sum construction from renormings of $\ell_1$}\\[2pt]
{\normalsize\itshape with further results on Banach lattices and dimension functions}\\[11pt]
{\large Ma Yicen}\\[8pt]
{\normalsize\bfseries AP--BAP conclusion obtained: October 30, 2026}\\[2pt]
{\small Revised and expanded: November 5, 2026}
\end{center}
\thispagestyle{plain}

\begin{abstract}
The author obtained the dual AP--BAP separation on October 30, 2026,
and recorded the conclusion and its proof in a dated five-page research
note. We give a direct construction of a separable Banach space $X$
such that $X$ and its continuous dual $X^*$ both have the approximation
property (AP) and both fail the bounded approximation property
(BAP). Take $X$ to be the $c_0$-sum of equivalent renormings of real
$\ell_1$ from OpenAI's quantitative renorming theorem. Finite coordinate
sums give AP, while their unbounded approximation constants rule out
BAP on both sides of the dual pairing. A complexification gives the
same separation over the complex field. The real predual has the weak
metric approximation property in every equivalent norm.
Further results concern strong generation in Bochner $L^1$ spaces,
disjoint projections in mixed norms, norm attainment in free Banach
lattices, and the simplex and density questions for dimension functions.
\end{abstract}

\section{The dual AP--BAP problem}\label{ap:intro}
If a continuous dual $X^*$ has the approximation property, can the
finite-rank approximations to its identity always be chosen with one
uniform norm bound? The problem concerns the given dual norm.
Figiel and Johnson separated AP from BAP for general Banach spaces
\cite{fj}. See Oja~\cite[Problem~2.1 and p.~221]{oja} for the
dual-space questions.

\paragraph{Original preprint.}
The author's five-page Zenodo preprint~\cite{ap:original}, dated
October~30, 2026, gives a negative answer to the dual AP--BAP
question conditional on OpenAI's quantitative renorming
result~\cite[Corollary~1.3]{oai}. The present paper is a revised and
expanded exposition dated November~5, 2026.

A space $Z$ has AP if, for each norm compact $K\subset Z$ and each
$\varepsilon>0$, a finite-rank bounded operator $T:Z\to Z$ satisfies
$\sup_{z\in K}\|Tz-z\|<\varepsilon$. It has $\lambda$-BAP if one
can always take $\|T\|\le\lambda$, BAP if some finite $\lambda$
suffices, and the metric approximation property (MAP) if
$\lambda=1$ suffices.

\begin{theorem}[AP without BAP on both sides of a dual pair]\label{ap:main}
For every integer $n\ge1$, let $E_n$ be an equivalent renorming of
real $\ell_1$ that fails $n$-BAP, as supplied by
\cite[Corollary~1.3]{oai}. Set
\[
 X=\left(\bigoplus_{n\ge1}E_n\right)_{c_0},
 \qquad X^*=\left(\bigoplus_{n\ge1}E_n^*\right)_{\ell_1}.
\]
Then $X$ is separable, $X^*$ is nonseparable, and each of $X$ and
$X^*$ has AP and fails BAP. A separable complex space with the same
two-sided separation is obtained in Corollary~\ref{ap:complex}.
\end{theorem}

The input gives a whole sequence of norms. For any one coordinate
the identity still admits uniformly bounded approximation. No single bound works for every coordinate, since
compressing an approximant for $X$ to $E_n$ preserves its bound.
In contrast, AP on a single compact set uses only finitely many
coordinates. Lemma~\ref{ap:dual} transfers the obstruction to
bounded approximation from $X$ to $X^*$.

\paragraph{Further results.}
Section~\ref{sec:sw} answers the strong-generation question of
Schl\"uchtermann and Wheeler~\cite{sworig,rodswcg,rodproj} by constructing
weakly compact diagonals in a Bochner space with Talagrand's Schur
dual as its range. Section~\ref{dj:section} identifies the diffuse
part of a mixed norm and constructs positive projections for sampling
kernels. Theorem~\ref{fr:main} characterizes scalar norm attainment in
free $p$-convex Banach lattices. Sections~\ref{bh:section} and
\ref{sx:sec} treat the two Blackadar--Handelman questions. A positive contraction separates a dimension function from all lower
semicontinuous ones, and the full dimension-function space, whenever
nonempty, is a compact Choquet simplex.

We write $B_Z$ for the closed unit ball of $Z$. Banach spaces and
lattices are real unless complex scalars are specified. Bochner
integrability includes strong measurability.

\section{A direct construction from renormings of \texorpdfstring{$\ell_1$}{l1}}\label{sec:ap}
OpenAI's quantitative renorming theorem~\cite[Corollary~1.3]{oai}
provides, for each $n\ge1$, a norm on real $\ell_1$ satisfying
\begin{equation}\label{ap:input}
 \tfrac12\|x\|_1\le\|x\|_{E_n}\le B_n\|x\|_1,\qquad
 B_n=300n^2+150n+5,
\end{equation}
for which $E_n$ fails $n$-BAP. Conjugating the coordinate truncations
of standard $\ell_1$ shows that $E_n$ still has $2B_n$-BAP.
The nonexistence of the bound $n$ is the substantive external input
to Theorem~\ref{ap:main}.

\subsection{Descent to the predual and the coordinate argument}
\begin{lemma}\label{ap:dual}
If $V^*$ has $\lambda$-BAP, then $V$ has $\lambda$-BAP.
\end{lemma}
\begin{proof}
Write $J:V\to V^{**}$ for the canonical embedding. Fix
$x_1,\ldots,x_k\in V$, finitely many test functionals in $V^*$,
and a tolerance. Choose a finite-rank $T:V^*\to V^*$ with
$\|T\|\le\lambda$ approximating the identity on those functionals.
By the principle of local reflexivity, for each $s>0$ there is
$L:T^*(V^{**})\to V$ with $\|L\|\le1+s$ that preserves evaluation by
all the given test functionals on this finite-dimensional domain.
Set $R=(1+s)^{-1}LT^*J$. Then $\|R\|\le\lambda$ and
\[
 h(Rx_i)=(1+s)^{-1}(Th)(x_i)
\]
for each test functional $h$. Taking the approximation error and
$s$ arbitrarily small puts $(x_i)_{i=1}^k$ in the weak closure in
$V^k$ of $\{(Rx_i)_{i=1}^k:R\text{ finite rank},\ \|R\|\le\lambda\}$.
This set is convex, so its norm closure is the same. Hence the
identity is approximable on finite sets with bound $\lambda$. A finite
net in a compact set and the bound $\|I-R\|\le1+\lambda$ give the
compact-set assertion.
\end{proof}

\begin{proof}[Proof of Theorem~\ref{ap:main}]
The dual identification is isometric: restriction to the summands
of $X$ defines a summable family of coordinate norms, as is seen by
testing on finitely supported vectors with separately chosen norming
signs and nearly norming coordinates. Conversely each such summable
family defines a functional by the coordinate pairing.

The spaces $E_n$ are isomorphic to $\ell_1$, and $E_n^*$ to
$\ell_\infty$. Both standard sequence spaces have AP. For
$\ell_1$ use the finite coordinate truncations. For $\ell_\infty$,
given finitely many bounded sequences, partition $\mathbb N$ into
finitely many cells on each of which every one of the sequences has
oscillation below the prescribed tolerance. Replacing a sequence on
each cell by its value at a representative gives a finite-rank
contraction approximating the chosen sequences. Uniform boundedness
passes this approximation to compact sets. AP is preserved by linear
isomorphisms, so every summand on each side has AP.

In either of the two sums the initial coordinate truncations are
contractions that converge uniformly on each norm compact $K$.
Choose $N$ with tail error $<\varepsilon/2$, then approximate
the identity on the compact image of $K$ in each of the first $N$
coordinates within $\varepsilon/(2N)$. These coordinate maps,
extended by zero on the tail, form a single finite-rank operator
with error $<\varepsilon$.

Every coordinate is a contractive retract of its ambient sum.
If $X$ had $\lambda$-BAP, composition with the coordinate injection
and projection would give $\lambda$-BAP on $E_n$ for all $n$. Taking
$n\ge\lambda$ contradicts the construction. Lemma~\ref{ap:dual}
then rules out BAP for $X^*$. Finally, $X$ is a countable $c_0$-sum of
separable spaces, while $X^*$ contains a coordinate copy of $\ell_\infty$.
\end{proof}

\begin{corollary}\label{ap:complex}
There is a separable complex Banach space $Z$ such that $Z$ and
its complex dual $D$ both have AP and fail BAP.
\end{corollary}
\begin{proof}
Complexify $X$ with
$\|x+iy\|_Z=\sup_{\theta\in\mathbb R}\|x\cos\theta-y\sin\theta\|_X$.
Its underlying real space is isomorphic to $X\oplus X$, and $X$ is
a contractive real retract via the real-part map. For $h\in X^*$
set $Ih(x+iy)=h(x)+ih(y)$. For $F\in D$ put $PF(x)=\operatorname{Re}F(x)$.
The maps $I,P$ are real-linear contractions with $PI$ the identity.
Writing $F=Ia+iIb$, with $a,b$ the real and imaginary parts of $F|_X$,
identifies the underlying real space of $D$ isomorphically with
$X^*\oplus X^*$. Both underlying real spaces therefore have AP.
To obtain complex AP on a compact set $K$ in either space, approximate
the identity on $K\cup iK$ by a real finite-rank map $A$, then take
$z\mapsto (Az-iA(iz))/2$. Conversely, bounded complex finite-rank
approximants would also be bounded real approximants. The two contractive
retractions would give BAP on $X$ or on $X^*$, respectively,
contradicting Theorem~\ref{ap:main}.
\end{proof}

\subsection{Weak metric approximation on the predual}
\begin{corollary}\label{ap:weak}
The real space $X$ in Theorem~\ref{ap:main} has the weak MAP in
every equivalent norm and the BAP in none.
\end{corollary}
\begin{proof}
The weak MAP tests finite-rank approximants after composition with a
weakly compact operator. More precisely, for each weakly compact
$T:X\to Z$ one requires a net $(S_i)$ of finite-rank operators on
$X$ converging to the identity uniformly on compact sets and satisfying
$\limsup_i\|TS_i\|\le\|T\|$. Lima and Oja~\cite[Theorem~4.2]{ap:weakmap}
characterize AP on $X^*$ by the weak MAP on
$X$ in every equivalent norm. Apply Theorem~\ref{ap:main} and use
the fact that existence of a finite BAP constant is invariant under
equivalent renorming.
\end{proof}

\section{Independent branches in a Bochner space}\label{sec:sw}
A space $Z$ is \emph{strongly weakly compactly generated} (SWCG)
if a weakly compact set $G\subset Z$ has the following property:
for each weakly compact $K\subset Z$ and $\varepsilon>0$, some
integer $n$ satisfies $K\subset nG+\varepsilon B_Z$. For strong
weak precompact generation (SWPG), both $G$ and the test sets $K$
are weakly precompact. A bounded set is weakly precompact if every
sequence in it has a weakly Cauchy subsequence.

\subsection{The exact cost of covering diagonals}
For a bounded sequence $v=(v_j)$ in a real Banach space set
\[
 \delta(v)=\sup_{z^*\in B_{Z^*}}
  \bigl(\limsup_j z^*(v_j)-\liminf_j z^*(v_j)\bigr).
\]
A bounded set is $\eta$-weakly precompact if every sequence in it
has a subsequence $v$ with $\delta(v)\le\eta$. Following
Rodr\'iguez~\cite{rodpre}, $Z$ has property $\mathfrak{KM}_w$ if
there are bounded sets $M_{n,p}$, $n,p\ge1$, such that
$M_{n,p}$ is $1/p$-weakly precompact and, for each weakly precompact
$K$ and each $p$, some $M_{n,p}$ contains $K$.
Every SWCG space is SWPG and every SWPG space has $\mathfrak{KM}_w$
\cite{rodpre}. For the latter implication take
$M_{n,p}=nG+(2p)^{-1}B_Z$ for a weakly precompact generator $G$.

Let $\mathfrak d$ be the least cardinality of an eventually
dominating family in $\mathbb N^{\mathbb N}$, with
$\mathbb N=\{1,2,\ldots\}$. For each $g$ at least one member
$f$ must satisfy $g(m)\le f(m)$ for all sufficiently large $m$. Closing
such a family under finite modifications makes it pointwise
dominating without increasing its cardinality. No countable family
$(f_i)$ eventually dominates $m\mapsto 1+\max_{i\le m} f_i(m)$.

\begin{theorem}\label{sw:array}
Let $(z_n^{(m)})_{m,n\ge1}$ be bounded in $Z$. Suppose every row
satisfies $\|\sum_n a_nz_n^{(m)}\|\ge c\sum_n|a_n|$ for finitely
supported scalars and a common $c>0$. Suppose also that
$z_{\phi(m)}^{(m)}\rightharpoonup0$ for each strictly increasing
$\phi:\mathbb N\to\mathbb N$. Set
$K_\phi=\{0\}\cup\{z_{\phi(m)}^{(m)}:m\ge1\}$.
\begin{enumerate}
\item For $0\le\eta<2c$, the least cardinality of a family of
$\eta$-weakly precompact sets such that each $K_\phi$ is contained
in one member is $\mathfrak d$.
\item For $0\le r<c$, the least cardinality of a family $\mathcal H$
of weakly compact sets such that each $K_\phi$ is contained in
$H+rB_Z$ for some $H\in\mathcal H$ is also $\mathfrak d$.
\end{enumerate}
Every Banach space containing $Z$ isomorphically fails $\mathfrak{KM}_w$.
\end{theorem}

\begin{proof}
On the span of any subsequence of any row, alternating values
$+c,-c,+c,\ldots$ prescribe a functional of norm at most one.
Hahn--Banach extends it to $Z$. Thus each row meets a given
$\eta$-weakly precompact set $M$ at finitely many indices if
$\eta<2c$. Put
$b_M(m)=\max(\{n:z_n^{(m)}\in M\}\cup\{1\})$.
If $K_\phi\subset M$, then $\phi\le b_M$ pointwise. Every
$g\in\mathbb N^{\mathbb N}$ has a strictly increasing majorant, so
a covering family has cardinality at least $\mathfrak d$.
The second assertion has the same lower bound: every $H+rB_Z$ is
$2r$-weakly precompact. Choose a weakly convergent subsequence of
the $H$-components and bound the scalar oscillation of the errors by $2r$.

For the upper bounds take a pointwise dominating family of size
$\mathfrak d$ with $f(m)\ge m$. Set
\[
 C_f=\{0\}\cup\{z_n^{(m)}:m\le n\le f(m)\}.
\]
Each row contributes finitely many elements. Any sequence in $C_f$
either has a constant subsequence or has a subsequence with increasing
row indices $m_j$. In the latter case we can also make the nonnegative
integers $n_j-m_j$ nondecreasing. Extend $\phi(m_j)=n_j$ by taking
$\phi(m)=m$ before $m_1$ and $\phi(m)-m=n_j-m_j$ on $[m_j,m_{j+1})$.
This is an increasing diagonal, so our subsequence converges weakly to
zero. Eberlein--\v Smulian gives the weak compactness of $C_f$, and
$\phi\le f$ implies $K_\phi\subset C_f$.

An embedding $J$ carries the array to one with the same diagonal
property and row lower bound $c/\|J^{-1}\|$, where the inverse is
taken on the range. Each $J(K_\phi)$ is weakly compact. Taking
$1/p<2c/\|J^{-1}\|$ rules out the countable family required by
$\mathfrak{KM}_w$.
\end{proof}

\subsection{The range space and the collision estimate}
Let $\mathcal T=\bigcup_{n\ge1}\{0,1\}^n$ be the dyadic tree and write
$\tau\succeq\sigma$ for extension of a word. Complete $c_{00}(\mathcal T)$
to a space $E$ under
\begin{equation}\label{sw:tree}
 \|b\|_E=\sup_{r\ge1}
 \left(\sum_{|\sigma|=r}\sup_{\tau\succeq\sigma}|b(\tau)|^2\right)^{1/2}.
\end{equation}
Talagrand's theorem gives that $X=E^*$ is separable and Schur
\cite{tal,roddelta}. In particular weakly compact subsets of $X$ are
norm compact. To see that $X$ is SWCG, take a dense sequence $(v_j)$
in its unit ball and the operator $S:\ell_2\to X$,
$Sa=\sum_j2^{-j}a_jv_j$. Its finite-rank truncations converge in
operator norm. Its range is dense and $S(B_{\ell_2})$ is compact,
since $S$ is compact and $B_{\ell_2}$ is weakly compact.
Approximating a finite net of a norm compact set by range vectors
shows that the integer dilates of $S(B_{\ell_2})$ strongly generate $X$.

Write $x_\sigma=e_\sigma^*$. These coordinates have norm one.
Along an infinite branch $\beta$, $(x_{\beta|n})$ is an isometric
$\ell_1$ basis: put the signs of any finite scalar sequence on the
corresponding branch nodes and zero elsewhere. The resulting test
vector has norm at most one in \eqref{sw:tree}. For words
$\tau_1,\ldots,\tau_M$ of length at least $r$, let $N_\sigma$ count
those with prefix $\sigma$ of length $r$. Then
\begin{equation}\label{sw:prefix}
 \left\|\sum_{k=1}^M x_{\tau_k}\right\|_X
 \le \left(\sum_{|\sigma|=r}N_\sigma^2\right)^{1/2}.
\end{equation}
Indeed evaluation on $b\in B_E$ bounds the sum by
$\sum_{|\sigma|=r} N_\sigma\sup_{\tau\succeq\sigma}|b(\tau)|$.
Apply Cauchy--Schwarz.

\begin{theorem}\label{sw:main}
For the space $X$ above, $Y=L^1([0,1],X)$ fails $\mathfrak{KM}_w$
and admits no isomorphic embedding into an SWPG space. In particular
SWCG does not pass from $X$ to $L^1([0,1],X)$.
\end{theorem}
\begin{proof}
Use independent fair bits on
$\Omega=\{0,1\}^{\mathbb N\times\mathbb N}$ with product probability.
In row $m$ put
$u_n^{(m)}(\omega)=x_{(\omega_{m,1},\ldots,\omega_{m,n})}$.
These are simple functions. The branch calculation gives
$\|\sum_n a_nu_n^{(m)}\|_{L^1}=\sum_n|a_n|$.
Consider an increasing diagonal $v_m=u_{\phi(m)}^{(m)}$ and a
subsequence $v_{m_j}$. For $M$ terms of lengths at least $r$, the
prefixes are independent uniform words. Counting equal ordered pairs gives
\[
 \mathbb E\sum_{|\sigma|=r}N_\sigma^2=M+M(M-1)2^{-r}.
\]
Integration of \eqref{sw:prefix} and Jensen's inequality bound the
norm of their sum by the square root of this expression. Put
$r=\lceil\log_2 N\rceil$ for $N\ge2$. Bound the first $r$ terms by
their norms. The remaining $M=N-r$ terms have length at least $r$,
since $\phi(m_j)\ge m_j\ge j$. Hence
\[
 \left\|N^{-1}\sum_{j=1}^N v_{m_j}\right\|_{L^1}
 \le \frac{\lceil\log_2 N\rceil}{N}+\sqrt{\frac2N}.
\]
This holds for every subsequence. A functional witnessing failure
of weak nullity would, after changing sign, be bounded below by a
positive constant on some subsequence. This contradicts the estimate.
Apply Theorem~\ref{sw:array} with $c=1$. Enumerating the coordinates
of $\Omega$ identifies product probability with Lebesgue measure
modulo null sets by binary expansion.
\end{proof}

Rodr\'iguez used the same range space to obtain a weakly compact
family failing uniform weakly compact control of its values
\cite{roddelta}. Here the independent copies allow selection against
any countable family of control sets. For this normalized array,
the covering numbers in Theorem~\ref{sw:array} drop to one at
$\eta=2$ and $r=1$: the unit ball gives the first cover and $\{0\}$
gives the second.

For a Banach space $Z$ let $\mathrm{SG}(Z)$ be the least size of a
family of weakly compact sets $H$ such that every weakly compact
$K\subset Z$ is contained, for each $\varepsilon>0$, in
$H+\varepsilon B_Z$ for some family member. The second part of
Theorem~\ref{sw:array}, applied to the image array, gives
$\mathrm{SG}(Z)\ge\mathfrak d$ in every space $Z$ containing an
isomorphic copy of $Y$. The matching upper bound in that theorem
concerns the diagonal family specifically.

\begin{corollary}\label{sw:atomic}
For the fixed space $X$ in Theorem~\ref{sw:main} and any probability
measure $\nu$, the following are equivalent: $\nu$ is purely atomic
modulo null sets, $L^1(\nu,X)$ is SWCG, $L^1(\nu,X)$ is SWPG, and
$L^1(\nu,X)$ has $\mathfrak{KM}_w$.
\end{corollary}
\begin{proof}
A purely atomic probability measure has at most countably many atoms
modulo null sets. Sending a Bochner function to its values on the
atoms multiplied by their masses identifies $L^1(\nu,X)$ with a
finite or countable $\ell_1$-sum of copies of $X$. Such sums preserve
SWCG~\cite[Proposition~2.9]{sworig}. Apply the implications at the
start of this section. Otherwise the atomless part $A$ has mass $a>0$.
Bisect its measurable subsets recursively into equal masses.
Under normalized measure on $A$, the resulting bits define a measurable
map $\psi:A\to\Omega$ with product distribution. The map sending $f$
to $a^{-1}f\circ\psi$ on $A$ and zero elsewhere is isometric on simple
functions and hence on $L^1(\Omega,X)$. Theorem~\ref{sw:array} now
excludes $\mathfrak{KM}_w$.
\end{proof}

\section{Diffuse rows and complemented disjoint sequences in mixed norms}\label{dj:section}

Let $(T,\tau)$ be a probability space, let $S$ be standard Borel, and let
$t\mapsto\nu_t$ be a probability kernel on $S$. Fix $1<q<p<\infty$ and
write $r=q/p$. We identify measurable functions modulo the null sets of the
barycentre $\mu=\int_T\nu_t\,d\tau(t)$ and use the lattice norm
\[
 \|f\|_E^q=\int_T\left(\int_S |f|^p\,d\nu_t\right)^r d\tau(t).
\]
The diffuse contribution can be recovered from the norm and order of $E$. In a purely atomic model the masses may depend on earlier observations, and a comparison with independent marks yields positive projections. These two results come from working notes
\cite{dj:diffnote,dj:adaptnote}. Cone embeddings and random-measure
representations in function lattices have been studied in~\cite{dj:hk}.

\subsection{Recovering the diffuse part}

Write $\nu_t=\nu_t^a+\nu_t^d$ for the atomic and diffuse parts. These are
kernels. Indeed, choose increasing finite Borel partitions $\mathcal D_j$
separating points of $S$. If $D_j(s)$ is the cell containing $s$, then
$\nu_t(\{s\})=\lim_j\nu_t(D_j(s))$ is jointly measurable, and restricting
$\nu_t$ to the zero set of this function gives $\nu_t^d$.

\begin{theorem}\label{dj:partition}
For $0\le h\in E$ put
\[
 L_N(h)=N^{r-1}\sup\left\{\sum_{i=1}^N\|h_i\|_E^q:
             h_i\ge0,\ \sum_{i=1}^N h_i\le h\right\}.
\]
Then $\displaystyle\lim_{N\to\infty}L_N(h)=
\int_T(\int_S h^p\,d\nu_t^d)^r\,d\tau(t)$.
The same limit results when the supremum is restricted to
$h_i=h\mathbf{1}_{A_i}$ for measurable partitions into at most $N$ cells,
padded with empty sets.
\end{theorem}

\begin{proof}
For any finite positive measure $\lambda$ and admissible $h_i$, concavity and
$\sum_i h_i^p\le h^p$ give
\[
 N^{r-1}\sum_i\left(\int h_i^p\,d\lambda\right)^r
 \le \left(\int h^p\,d\lambda\right)^r.
\]
Choose a family for each $N$ within $1/N$ of the normalized supremum. The
inequalities for these countably many families hold rowwise outside a common
$\tau$-null set. Fix such a row with finite $p$-integral. On a finite set
$F$ of its atoms, subadditivity of $x^r$ bounds the unnormalized sum by
$\sum_{s\in F}\nu_t(\{s\})^r\sum_i h_i(s)^q
\le\sum_{s\in F}\nu_t(\{s\})^r h(s)^q$.
After normalization this tends to zero. The normalized contribution of the
remaining atoms is at most $(\int_{S\setminus F}h^p\,d\nu_t^a)^r$, which
can be made arbitrarily small. Thus the atomic contribution tends to zero
rowwise. The displayed bound supplies an integrable majorant. Subadditivity
for the two parts of each row and dominated convergence give the upper bound.

For the lower bound fix $N$. Assign the cells of $\mathcal D_j$ independently
and uniformly to $N$ colours, and let $A_i$ be the union of the cells of
colour $i$. For $\lambda_t=h^p\nu_t^d$ the random variable $\lambda_t(A_i)$
has mean $\lambda_t(S)/N$ and variance at most
$\sum_{D\in\mathcal D_j}\lambda_t(D)^2$. This sum tends to zero, since the
sets of pairs lying in the same cell decrease to the diagonal, which has
$\lambda_t\otimes\lambda_t$-measure zero. For almost every row the measures
are finite, so bounded convergence in probability gives convergence of the
expected $r$th powers. The displayed bound permits integration over the rows.
The expected normalized sum for the diffuse rows therefore tends to
$\int_T\lambda_t(S)^r\,d\tau(t)$. The supremum over deterministic partitions
for the full rows dominates each expectation, giving the claimed lower bound.
\end{proof}

It follows that if almost every row is purely atomic, there is no positive
linear $U:L^p[0,1]\to E$ with $\|Ug\|_E\ge c\|g\|_p$ for all $g\ge0$
and some $c>0$. Indeed, for the $N$ equal intervals $I_i$, take
$h=U\mathbf{1}$ and $h_i=U\mathbf{1}_{I_i}$. Then $L_N(h)\ge c^q$ for
every $N$, whereas Theorem~\ref{dj:partition} gives limit zero.

\subsection{Predictable atomic sampling}

Let $Y_i$ be $\mathcal F_i$-measurable marks in a standard atomless
probability space $(S,\lambda)$, each with conditional law $\lambda$ given
$\mathcal F_{i-1}$. Let $b_i\ge0$ be $\mathcal F_{i-1}$-measurable with
$\sum_i b_i=1$ almost surely. Use the kernel
$\nu_\omega=\sum_i b_i(\omega)\delta_{Y_i(\omega)}$ in the definition of $E$.
Its barycentre is $\lambda$ by conditioning. Let $(V_i)$ be independent
marks of law $\lambda$, independent of the entire original process. Set
\[
 \|f\|_R^q=\mathbb E\left(\sum_i b_i |f(V_i)|^p\right)^r.
\]
This is a rearrangement-invariant norm on its natural domain by the
mixed-norm triangle inequality and independence of the new marks. The
comparison below is related to estimates for tangent sequences
\cite{dj:tangent}. We give an elementary proof for these sampling kernels.

\begin{lemma}\label{dj:stop}
For adapted $0\le Z_i\le1$ and $a_i=\mathbb E(Z_i\mid\mathcal F_{i-1})$,
write $S_n=\sum_{i\le n}Z_i$ and $A_n=\sum_{i\le n}a_i$. Then
\[
 \tfrac12\mathbb E(1\wedge A_n)\le \mathbb E(1\wedge S_n)
 \le 2\mathbb E(1\wedge A_n).
\]
\end{lemma}
\begin{proof}
Let $A_0=0$. Choose predictable $c_i=\min\{1,(1-A_{i-1})_+/a_i\}$ if
$a_i>0$ and $c_i=1$ otherwise. Then $\sum_{i\le n}c_i a_i=1\wedge A_n$
and $1\wedge S_n\le\mathbf1_{\{A_n\ge1\}}+\sum_{i\le n}c_iZ_i$.
Taking expectations gives the upper bound.
For the reverse inequality adjoin independent uniform $[0,1]$ random
variables $U_i$ independent of the original process. Declare trial $i$ a
success if $U_i\le Z_i$, and let $H$ be the event of any success by time $n$.
Conditioning on the original process gives
$\mathbb P(H)=\mathbb E[1-\prod_{i\le n}(1-Z_i)]\le\mathbb E(1\wedge S_n)$.
Writing $G_i$ for no success before $i$, conditioning on the enlarged past
gives $\mathbb E\sum_{i\le n}\mathbf1_{G_i}a_i=\mathbb P(H)$ by summing
the disjoint first-success probabilities. Since
$1\wedge A_n\le\mathbf1_H+\sum_{i\le n}\mathbf1_{G_i}a_i$, the lower bound follows.
\end{proof}

\begin{theorem}\label{dj:predict}
The norms satisfy $4^{-1}\|f\|_R^q\le\|f\|_E^q\le4\|f\|_R^q$.
For every countable measurable partition $(C_j)$ with positive cell measures,
conditional expectation onto the partition defines a positive projection on
$E$ of norm at most $16^{1/q}$. Its range is the closed span of the cell
indicators. In particular $E$ has a complemented disjoint sequence.
\end{theorem}
\begin{proof}
Fix $t>0$. Apply Lemma~\ref{dj:stop} to
$Z_i=1\wedge(tb_i|f(Y_i)|^p)$. Its predictable mean is
$a_i=\int_S [1\wedge(tb_i|f(x)|^p)]\,d\lambda(x)$.
Conditional on all the masses, the independent marks $V_i$ have the same
formula for the predictable means, with the masses fixed. Applying the lemma
in this conditional model and integrating yields the same two bounds in
terms of $\mathbb E(1\wedge\sum_{i\le n}a_i)$. Hence the two truncated
expectations, for the original and independent marks, are within a factor
of four. Taking the minimum of one and a sum commutes with truncating each
nonnegative summand at one. Letting $n\to\infty$ and using Tonelli with
\[
 x^r=r(1-r)\int_0^\infty(1\wedge tx)t^{-1-r}\,dt
\]
gives the norm comparison.

Write $P f$ for the function equal on $C_j$ to
$\lambda(C_j)^{-1}\int_{C_j} f\,d\lambda$. It is defined for $f\in E$
since $\|f\|_{L^q(\lambda)}\le\|f\|_E$ by Jensen for the probability rows.
Condition on the masses and the cell labels of $V_1,\ldots,V_n$. The
conditional expectation of the vector $(b_i^{1/p}f(V_i))_{i\le n}$ is
$(b_i^{1/p}P f(V_i))_{i\le n}$. Jensen for the convex function
$v\mapsto\|v\|_{\ell_p^n}^q$, followed by monotone convergence, gives
$\|P f\|_R\le\|f\|_R$. The comparison now gives the asserted bound in
$E$. The averaging formula gives positivity and $P^2=P$. For a function
in the range, truncating to finitely many cells converges in $E$ by
dominated convergence first within and then across the rows. Thus the range
is the stated closed span. Choose an infinite partition, which exists by
atomlessness, and normalize its indicators to obtain a disjoint basic
sequence (coordinate truncations are contractive in the lattice norm).
\end{proof}

The general question asking whether every infinite-dimensional separable Banach lattice
admits a disjoint sequence with complemented span is recorded in
\cite[Question~8]{dj:survey}. The preceding projection applies to the stated
predictable-kernel model.

\subsection{Power weights determine an intermediate exponent}

The following computation from~\cite[p.~5]{dj:rajnote} identifies the
range geometry for deterministic power weights. In the translated model it
also describes exactly which convolution expression controls the norm.

\begin{theorem}\label{dj:power}
Fix $q<s<p$ and $b_i=c i^{-p/s}$, where $c$ normalizes the sum to one.
For independent marks of common law $\kappa$, put
$F(f)=(\mathbb E(\sum_i b_i|f(Y_i)|^p)^{q/p})^{1/q}$.
There are constants depending only on $p,q,s$ such that
$F(f)\asymp\|f\|_{L^s(\kappa)}$, uniformly over probability laws $\kappa$.
On the circle let $T$ have Haar law $m$, and conditionally on $T=t$ let the
marks be independent with law $\nu(\cdot-t)$ for an arbitrary probability
$\nu$. Writing $K_\nu g(t)=\int g(t+y)\,d\nu(y)$, the resulting sampling
norm satisfies
\[
 \|f\|_E\asymp \left(\int (K_\nu |f|^s)^{q/s}\,dm\right)^{1/q}.
\]
For $\nu=m$ this is an equivalent norm on $L^s(m)$, and every infinite
measurable partition into positive-measure cells gives a positive,
contractively complemented disjoint copy of $\ell_s$ in the sampling norm.
\end{theorem}
\begin{proof}
Set $\alpha=s/p$, $r=q/p$, and first normalize $\int|f|^s\,d\kappa=1$.
The function $x\mapsto1-\exp(-cvx^{-1/\alpha})$ is decreasing on
$(0,\infty)$ with integral $A v^\alpha$, where
$A=c^\alpha\Gamma(1-\alpha)$ by substitution and integration by parts.
Comparison with the sums over unit intervals therefore gives
\[
 A v^\alpha-1\le B(v):=\sum_{i\ge1}(1-e^{-vb_i})\le A v^\alpha.
\]
Put $X=\sum_i b_i|Y_i'|^p$, where $Y_i'=f(Y_i)$, and
$h_i(v)=\int (1-e^{-vb_i|f|^p})\,d\kappa$. Integrating the bounds for
$B(v|f|^p)$ shows that $Av^\alpha-1\le\sum_i h_i(v)\le Av^\alpha$.
By independence $\mathbb E e^{-vX}=\prod_i(1-h_i(v))$. Hence
$1-\mathbb E e^{-vX}\le\min\{1,Av^\alpha\}$, whereas for
$v\ge v_0=(2/A)^{1/\alpha}$ it is at least $1-e^{-1}$. In the identity
\[
 \mathbb E X^r=\frac{r}{\Gamma(1-r)}\int_0^\infty
       (1-\mathbb E e^{-vX})v^{-1-r}\,dv
\]
these bounds give finite upper and strictly positive lower constants since
$0<r<\alpha$. Homogeneity and monotone truncation establish the comparison
for all functions. Conditional application to each translated law, followed
by integration in $t$, gives the formula for $E$. If $\nu=m$, the row law
is constant. Averaging on any partition is a contraction in the sampling
norm by the conditional Jensen argument in Theorem~\ref{dj:predict}, since
the weights are deterministic and the marks independent. Its range, given
the equivalent $L^s$ norm, is isomorphic to $\ell_s$ by normalizing the cell
indicators in $L^s(m)$.
\end{proof}

\section[Scalar norm attainment in the free p-convex Banach lattice]{Scalar norm attainment in the free $p$-convex Banach lattice}\label{fr:section}
Let $E$ be a real Banach space and $1\le p<\infty$. Write $F_p(E)=\mathrm{FBL}^{(p)}[E]$ for the free $p$-convex Banach lattice with convexity constant one. In its concrete realization the generators are the functions $\delta_x(h)=h(x)$ on $E^*$, for $x\in E$. The vector sublattice generated by these functions is dense in $F_p(E)$. The norm formula of~\cite{fr:jlt} is
\begin{equation}\label{fr:norm}
 \|f\|_{F_p} =\sup\left\{\left(\sum_{j=1}^m|f(h_j)|^p\right)^{1/p} :
  \sup_{x\in B_E}\sum_{j=1}^m |h_j(x)|^p\le 1,\quad h_j\in E^*,\ m\ge1\right\}.
\end{equation}
Every such $f$ is positively homogeneous. For $u\in E^*$ define $\widehat u(f)=f(u)$. This is the lattice homomorphism extending $u$ and $\|\widehat u\|=\|u\|$. Indeed, a one-term test in \eqref{fr:norm} gives the upper bound, and testing on $\delta_x$ gives the reverse inequality since $\|\delta_x\|_{F_p}=\|x\|$. We use norm attainment on the closed unit ball, which includes the zero case.

\subsection{Scalar evaluations and AM-space targets}
\begin{theorem}\label{fr:main}
For every real Banach space $E$, every $1\le p<\infty$ and every $u\in E^*$, the functional $\widehat u$ attains its norm on $F_p(E)$ if and only if $u$ attains its norm on $E$.
\end{theorem}
This is the equivalence asked for in~\cite[Question~6.3]{fr:bil}. At $p=1$ it gives~\cite[Conjecture~5.5]{fr:dan}, whose authors established the assertion for several classes of spaces. The bidual of $E$ supplies a norming point for $u$ without any attainment assumption. The point that requires proof here is that equality at a free lattice element gives a norming point in the canonical image of $E$ itself.

\begin{lemma}[Exact domination, {\cite[Proposition~7.7]{fr:jlt}}]\label{fr:dom}
Let $f\in B_{F_p(E)}$ and equip $K=B_{E^{**}}$ with $\sigma(E^{**},E^*)$. There is a regular Borel probability measure $\mu$ on $K$ such that
\begin{equation}\label{fr:ineq}
  |f(h)|^p\le\int_K|z(h)|^p\,d\mu(z)\qquad(h\in E^*).
\end{equation}
\end{lemma}
\begin{proof}
The probability measures on the compact Hausdorff space $K$ form a weak-star compact convex subset of $C(K)^*$. Each constraint in \eqref{fr:ineq} defines a closed subset of this set. Fix $h_1,\ldots,h_m$ and let $D\subset\mathbb R^m$ be the compact convex set of vectors
\[
\left(\int_K |z(h_j)|^p\,d\nu(z)-|f(h_j)|^p\right)_{j=1}^m
\]
as $\nu$ ranges over these measures. If $D$ were disjoint from the nonnegative orthant, strict separation would give $\lambda_j\ge0$ such that $\sup_{d\in D}\sum_j\lambda_j d_j<0$. Apply \eqref{fr:norm} to the scaled list $\lambda_j^{1/p}h_j$ and use positive homogeneity to get
\[
\sum_j\lambda_j |f(h_j)|^p\le\sup_{x\in B_E}\sum_j\lambda_j|h_j(x)|^p.
\]
The point masses at the canonical images of $x\in B_E$ now give a contradiction. Thus all finite collections of constraints hold simultaneously, and compactness yields $\mu$.
\end{proof}

\begin{proof}[Proof of Theorem~\ref{fr:main}]
If $u$ is normed by $x\in B_E$, then $\delta_x$ norms $\widehat u$. Conversely we may assume $\|u\|=1$ and choose $f\in B_{F_p(E)}$ with $|f(u)|=1$. Take $\mu$ as in Lemma~\ref{fr:dom}. Since $|z(u)|\le1$ on $K$, equality in \eqref{fr:ineq} at $u$ implies $|z(u)|=1$ for $\mu$-almost every $z$. Define $b\in E^{**}$ by
\[
 b(h)=\int_K z(u)z(h)\,d\mu(z).
\]
We have $\|b\|\le1$ and $b(u)=1$. Let $J_E:E\to E^{**}$ be the canonical embedding. If $b\notin J_E(E)$, the Hahn--Banach theorem on $E^{**}/J_E(E)$ gives $\xi\in B_{E^{***}}$ with $\xi|_{J_E(E)}=0$ and $\xi(b)=c>0$. Goldstine's theorem applied to $E^*$ gives a net $(h_i)\subset B_{E^*}$ with $z(h_i)\to\xi(z)$ for every $z\in E^{**}$. In particular $h_i\to0$ in $\sigma(E^*,E)$ and $b(h_i)\to c$.

Every element of $F_p(E)$ is weak-star continuous on bounded dual balls. In fact lattice polynomials in the $\delta_x$ have this property, and the one-term bound from \eqref{fr:norm} gives uniform convergence on each such ball under approximation in free norm. Fix $0<s<1/2$. For almost every $z$, real scalar arithmetic gives
\[
 |z(u-sh_i)|=1-sz(u)z(h_i).
\]
Taylor's estimate $(1-t)^p\le1-pt+C_p t^2$ for $|t|\le1/2$ holds with a finite constant independent of $i$ and $s$, and with $C_1=0$. Apply \eqref{fr:ineq}, integrate the estimate, and take the limit along the net. Weak-star continuity on $2B_{E^*}$ gives
\[
 1=|f(u)|^p=\lim_i |f(u-sh_i)|^p\le 1-psc+C_p s^2.
\]
After division by $s$, letting $s\downarrow0$ contradicts $c>0$. Thus $b=J_Ex$ for some $x\in B_E$, and $u(x)=b(u)=1$.
\end{proof}

\begin{corollary}[Operators with AM-space range]\label{fr:am}
Let $Y$ be a real AM-space, meaning that $\|v\vee w\|=\max(\|v\|,\|w\|)$ for $v,w\ge0$. For every bounded $T:E\to Y$, its canonical lattice homomorphism $\widehat T:F_p(E)\to Y$ attains its norm if and only if $T$ does.
\end{corollary}
\begin{proof}
The Kakutani representation gives an isometric lattice embedding $j:Y\to C(L)$ for a compact Hausdorff space $L$. The free extension has $\|\widehat T\|=\|T\|$ by the universal property, since an AM-space is $p$-convex with constant one. If $x\in B_E$ norms $T$, then $\delta_x$ norms $\widehat T$. Conversely suppose that $\|\widehat T f\|=\|T\|>0$ for $f\in B_{F_p(E)}$. Choose $t\in L$ with $|(j\widehat T f)(t)|=\|T\|$ and set $u(x)=(jTx)(t)$. By uniqueness of free extensions, $\widehat u(g)=(j\widehat T g)(t)$ for every $g$. Hence $\|u\|=\|T\|$ and $f$ norms $\widehat u$. Theorem~\ref{fr:main} gives $x\in B_E$ with $|u(x)|=\|T\|$, whence $\|Tx\|=\|T\|$.
\end{proof}

\subsection[An exact test on an L1 domain]{An exact test on an $L^1$ domain}\label{fr:boundary}
Let $(\Omega,\mathcal A,\mu)$ be a probability space and $m\ge1$. Each bounded operator $T:L^1(\mu)\to\ell_p^m$ has a coordinate representation
\[
 Tx=\int_\Omega x(t)g(t)\,d\mu(t),\qquad
 g=(g_1,\ldots,g_m),\quad g_j\in L^\infty(\mu).
\]
Let $p'$ be the conjugate exponent, with $p'=\infty$ at $p=1$. Duality gives $\|T\|=\sup_{a\in B_{\ell_{p'}^m}}\|\sum_j a_jg_j\|_\infty=\operatorname*{ess\,sup}_t\|g(t)\|_p$. The last identity follows by taking a countable dense subset of the finite-dimensional coefficient ball. The extension has coordinates $\widehat T f=(f(g_j))_{j=1}^m$ by the defining property of the scalar extensions.

\begin{proposition}\label{fr:kernel}
Suppose $\|T\|=1$ and set $A=\{t:\|g(t)\|_p=1\}$. Then $\widehat T$ attains its norm if and only if $\mu(A)>0$. The operator $T$ attains its norm if and only if some $a\in S_{\ell_{p'}^m}$ satisfies $\mu\{t:|a(g(t))|=1\}>0$.
\end{proposition}
\begin{proof}
For a list $(h_i)_{i=1}^k$ in $L^\infty(\mu)$, the admissibility constraint in \eqref{fr:norm} means exactly $\sum_i |h_i(t)|^p\le1$ almost everywhere, by the same norm identity applied with $g=(h_i)_{i=1}^k$.

First take $f\in B_{F_p(L^1(\mu))}$ with $\sum_j|f(g_j)|^p=1$. For $n\ge2$ set $D_n=\{t:\|g(t)\|_p\le1-1/n\}$. Append all $n^{-1}\mathbf1_{D_n}g_j$ to the list $(g_j)$. The enlarged list is admissible since $(1+n^{-p})^{1/p}(1-1/n)\le(1+1/n)(1-1/n)\le1$. The initial list already exhausts the unit bound in \eqref{fr:norm}, so $f(\mathbf1_{D_n}g_j)=0$ for all $j$. If $\mu(A)=0$, then $\mathbf1_{D_n}g_j\to g_j$ weak-star on a bounded ball by dominated convergence against $L^1(\mu)$. Bounded-ball weak-star continuity of $f$ gives $f(g_j)=0$ for every $j$, a contradiction.

Conversely let $\mu(A)>0$, set $\rho=\mathbf1_A/\mu(A)$ and $d\nu=\rho\,d\mu$. For each $j$ choose $z_j\in B_{L^{p'}(\nu)}$ with $\int z_j g_j\,d\nu=\|g_j\|_{L^p(\nu)}$, using the sign function when $p=1$. Extend $x_j=\rho z_j$ by zero off $A$, so $x_j\in L^1(\mu)$. The lattice polynomial $f=\bigvee_{j=1}^m |\delta_{x_j}|$ satisfies $f(h)\le\|h\|_{L^p(\nu)}$ by H\"older's inequality, with equality at every $g_j$. Integration of the pointwise admissibility constraint gives $\|f\|_{F_p}\le1$. Furthermore,
\[
 \|\widehat T f\|_p^p=\sum_j\|g_j\|_{L^p(\nu)}^p
  =\int_A\|g(t)\|_p^p\,d\nu(t)=1.
\]
For the final assertion, if $x\in B_{L^1(\mu)}$ norms $T$, choose a norming functional $a$ for $Tx$. Equality in $1=a(Tx)\le\int |x(t)|\,|a(g(t))|\,d\mu(t)\le1$ forces $|a(g(t))|=1$ for $|x|\,d\mu$-almost every $t$, hence on a set of positive $\mu$-measure. Conversely, the normalized signed indicator of that set norms $T$, taking its sign to be that of $a(g)$.
\end{proof}

\begin{corollary}\label{fr:arc}
For every $1<p<\infty$ there is a rank-two operator $T:L^1[0,1]\to\ell_p^2$ that fails to attain its norm although $\widehat T$ attains it at a lattice polynomial. For $p=1$, every $T:L^1(\mu)\to\ell_1^m$ attains its norm whenever its extension does.
\end{corollary}
\begin{proof}
For $p>1$ define $v(t)=(\cos t,\sin t)$ on $[0,1]$, and take $g(t)=v(t)/\|v(t)\|_p$. Its norm is identically one. The arc spans the plane, and continuity shows that no nonzero functional annihilates $g$ almost everywhere. Hence $T$ has rank two. Strict convexity of $\ell_p^2$ gives at most one antipodal pair at which a given unit functional has absolute value one. Each such pair meets this injective arc at most once. Proposition~\ref{fr:kernel} gives the two claims about attainment. For $p=2$, this kernel and the nonattainment of its integral operator are due to Kadets, Mart\'in, Mer\'i and Quero~\cite[Example~3]{fr:kmmq}. The proposition supplies attainment of the free extension.

When $p=1$, the finitely many sign vectors in $\{-1,1\}^m$ norm $\ell_1^m$. If $\mu(A)>0$, one of these functionals norms $g(t)$ on a subset of $A$ of positive measure. Apply the second assertion of Proposition~\ref{fr:kernel}, with the zero operator handled directly.
\end{proof}

\section{Failure of density for lower semicontinuous dimension functions}\label{bh:section}

Blackadar and Handelman asked whether lower semicontinuous dimension functions
are pointwise dense among all normalized dimension functions of a unital
complex $C^*$-algebra~\cite[following Theorem I.2.4]{bh:original}.
Finite radius of comparison is a sufficient condition for density
\cite[Theorem~8.9]{bh:aprt}. We construct a separable example in which a
single positive element separates the two sets. The bundle used in each
coordinate occurs in~\cite[Example~7.5]{bh:rr}. The choice of a sequence of
such coordinates and the separating element are taken from~\cite[Section~7]{bh:source}.

For a unital complex $C^*$-algebra $A$, write $W(A)$ for the ordered
abelian monoid of Cuntz classes of positive elements of finite matrix algebras
over $A$. Thus $[b]\le[c]$ means $v_j^*cv_j\to b$ in norm for a sequence
of compatible rectangular matrices, with zero padding allowed. Addition is
block sum. Let $u=[1_A]$. A normalized dimension function is an additive
order-preserving map $d:W(A)\to[0,\infty)$ with $d(u)=1$.
We write $d(b)$ for $d([b])$ and give $\mathrm{DF}(A)$ the topology of
pointwise convergence on $W(A)$. Its subset $\mathrm{LDF}(A)$ consists of
the dimension functions lower semicontinuous on positive matrices in their
norm topology. General states on $W(A)$ need not preserve suprema of
increasing sequences in the completed Cuntz semigroup $\mathrm{Cu}(A)$.

\begin{theorem}\label{bh:main}
There are a separable, unital, nuclear, stably finite $2$-subhomogeneous
$C^*$-algebra $A$ and a positive contraction $a\in A$ such that
\[
 \max_{d\in\mathrm{DF}(A)}d(a)=1,\qquad
 \max_{e\in\mathrm{LDF}(A)}e(a)=\tfrac12 .
\]
Thus lower semicontinuous dimension functions need not be pointwise dense.
\end{theorem}

\begin{proof}
For $n\ge1$ let $L_n$ be the tautological complex line bundle on
$\mathbb{CP}^n$, set $E_n=\mathbf{1}\oplus L_n$, and let
$B_n=\Gamma(\operatorname{End}E_n)$, the algebra of continuous bundle
endomorphisms. The symbol $\mathbf{1}$ denotes the trivial complex line.
Choose Hermitian metrics and let $p_n\in B_n$ project onto $\mathbf{1}$.
Take the minimal unitization of the $c_0$-sum and put
\[
 I=\bigoplus_{n\ge1}^{c_0} B_n,\qquad A=I^{\sim},\qquad
 \varepsilon_n=(n+1)^{-1},\qquad a=(\varepsilon_n p_n)_{n\ge1}.
\]
The decay of $\varepsilon_n$ puts $a$ in $I$. Each $B_n$ is a separable
homogeneous algebra with matrix fibers $M_2(\mathbb C)$. In an irreducible
representation of $A$, the central units of the summands act by zeros and
ones. If they all act by zero, the representation factors through $A/I
\cong\mathbb C$. Otherwise it factors through one $B_n$ and is a fiber
evaluation. Thus $A$ is $2$-subhomogeneous and type I, hence nuclear.
It is separable. The coordinate evaluations and quotient character separate
points at every matrix level, so every isometry at every matrix level is
unitary. This gives stable finiteness.

We first obstruct comparison with the order unit. Suppose
$(m+1)[a]\le m u$ for some integer $m\ge1$. Evaluation at coordinate
$n=m$ is a unital $*$-homomorphism, and a nonzero scalar multiple of a
projection is Cuntz equivalent to that projection. Hence
$(m+1)[p_m]\le m[1_{B_m}]$. For projections Cuntz subequivalence is
Murray--von Neumann subequivalence. The resulting partial isometry gives
a bundle embedding of $\mathbf{1}^{\oplus(m+1)}$ into
$E_m^{\oplus m}$ over $\mathbb{CP}^m$. Its orthogonal complement $G$
would have rank $m-1$. With $x=c_1(L_m)$, the Whitney product formula gives
\[
 c(G)=c(E_m^{\oplus m})=(1+x)^m .
\]
The degree-$2m$ term $x^m$ is nonzero in
$H^*(\mathbb{CP}^m;\mathbb Z)=\mathbb Z[x]/(x^{m+1})$.
Yet $c_m(G)=0$ by the rank of $G$. Thus $(m+1)[a]\nleq m u$ for every $m$.

We now use the ordered-monoid state criterion
\cite[Proposition~5.2.13]{bh:apt}. In a positively ordered abelian monoid
with order unit $u$, an element $y\le u$ satisfies $(k+1)y\le k u$ for
some $k\ge1$ if and only if every extended additive order-preserving
state $f$ normalized at $u$ satisfies $f(y)<1$. The order-unit condition
makes all such states finite on the monoid: $z\le Nu$ gives $f(z)\le N$.
It holds in $W(A)$, since a positive $N$-by-$N$ matrix is Cuntz
subequivalent to $1_{M_N(A)}$. In particular $[a]\le u$. The comparisons
we just excluded therefore supply $d\in\mathrm{DF}(A)$ with $d(a)=1$.
No dimension function can take a larger value on $a$.

Finally, by~\cite[Theorems~II.2.2, II.3.1 and II.4.9]{bh:original},
each $e\in\mathrm{LDF}(A)$ is $d_\tau$ for a tracial state $\tau$,
since $A$ is type I, and $d_\tau(a)=\lim_{k\to\infty}\tau(a^{1/k})$.
Let $e_n$ be the central unit of $B_n$ in $A$, and put
$w_n=\tau(e_n)$. Then $\sum_n w_n\le1$.
To see directly that $\tau(p_n)=w_n/2$, choose local unit sections of
$L_n$ and a finite partition of unity subordinate to their trivializing
cover. Multiply each section by the square root of its partition function
and extend by zero. The resulting bundle maps $v_i$ from the trivial
summand to $L_n$, taken to be zero on $L_n$, satisfy
$\sum_i v_i^*v_i=p_n$ and $\sum_i v_i v_i^*=e_n-p_n$.
Taking traces gives the assertion. For fixed $k$ the expansion of $a^{1/k}$
in its summands converges in norm. Monotone convergence for the resulting
nonnegative series yields
\[
 d_\tau(a)=\lim_{k\to\infty}\frac12\sum_{n\ge1}w_n\varepsilon_n^{1/k}
           =\frac12\sum_{n\ge1}w_n\le\frac12 .
\]
The normalized matrix trace after evaluation in any one fiber attains the
upper bound. The pointwise closed set $\{f\in\mathrm{DF}(A): f(a)\le1/2\}$
therefore contains $\mathrm{LDF}(A)$ and omits $d$.
\end{proof}

\section{The simplex of dimension functions}\label{sx:sec}
We keep the pointwise topology on $\mathrm{DF}(A)$ from Section~\ref{bh:section}. The simplex conjecture of Blackadar and Handelman~\cite{bh:original} concerns this full state space. Antoine, Perera, Robert and Thiel proved the stable-rank-one case~\cite{sx:sr}. Gardella and Perera give a later account in~\cite[Section~12]{sx:gp}.

\begin{theorem}\label{sx:main}
For any unital complex $C^*$-algebra $A$, the space $\mathrm{DF}(A)$, if nonempty, is a compact Choquet simplex. In particular the conclusion holds for every nonzero unital stably finite $A$.
\end{theorem}
We give the rank-space and Hilbert-module arguments of~\cite[Sections~2--6]{bh:source}. We use the compact-base characterization of a Choquet simplex: the real span of its positive cone, ordered by that cone, is a vector lattice. This characterization applies without metrizability.

\subsection{Additive functions and Sylvester ranks}
\begin{lemma}\label{sx:length}
Suppose $\mathcal C$ is an essentially small abelian category. The finite, nonnegative real functions $L$ on isomorphism classes that are additive on short exact sequences form a cone whose ordered span is a vector lattice. If every object is a subquotient of a finite sum of copies of $U$, the nonempty normalized base $L(U)=1$ is a compact Choquet simplex.
\end{lemma}
\begin{proof}
For two positive functions $L_1,L_2$, take
\[
 h(M)=\inf\sum_i\min\{L_1(Q_i),L_2(Q_i)\},
\]
where the infimum ranges over all finite filtrations of $M$ by subobjects, with successive quotients $Q_i$. For $0\to X\to M\to Y\to0$, intersect a filtration with $X$ and take its image in $Y$. The two quotients at each step fit into an exact sequence with middle term $Q_i$. Since $\min(a+c,b+d)\ge\min(a,b)+\min(c,d)$, each filtration costs at least $h(X)+h(Y)$. Concatenating a filtration in $X$ with a lifted filtration from $Y$ gives the reverse inequality. Thus $h$ is additive. The single-step filtration gives $h\le L_1,L_2$. For any additive lower bound $v$ in the span, additivity along each filtration gives $v(M)\le\sum_i\min\{L_1(Q_i),L_2(Q_i)\}$ even if $v$ has negative values. Thus $h$ is the meet. Translating any pair of elements in the span into the cone and translating its meet back gives a meet for the pair. Negation gives joins.

If $M$ is a subquotient of $U^{\oplus k_M}$, then $0\le L(M)\le k_M L(U)$. Hence $L(U)=0$ implies $L=0$, and normalization yields a base. On a skeleton it is closed in $\prod_M[0,k_M]$ by the additive identities. Compactness and the vector-lattice characterization give the assertion.
\end{proof}

For a unital ring $R$, a normalized Sylvester matrix rank $r$ assigns a finite nonnegative value to each rectangular matrix. It vanishes on zero, takes value one on $1_R$, is additive on block sums and satisfies $r(xy)\le\min(r(x),r(y))$ and
$ r\bigl(\begin{smallmatrix} x&z\\0&y\end{smallmatrix}\bigr)\ge r(x)+r(y)$.
Write $P(R)$ for its space with the pointwise topology.

\begin{proposition}\label{sx:rank}
If $P(R)$ is nonempty, it is a compact Choquet simplex.
\end{proposition}
\begin{proof}
Let $\mathcal D$ be a skeleton of finitely presented left $R$-modules and $H_X=\mathrm{Hom}_R(X,-)$. By the Malcolmson matrix/module correspondence~\cite{sx:hl}, ranks correspond affinely to normalized nonnegative functions $\rho$ on $\mathcal D$, additive on direct sums, satisfying
\[
 \rho(Z)\le\rho(Y)\le\rho(X)+\rho(Z)\quad
 \text{whenever }X\to Y\to Z\to0\text{ is exact}.
\]
A presentation on $m$ generators with relation matrix $x$ gives $\rho(M)=m-r(x)$. The additive-functor description of such ranks is also recorded by Virili~\cite{sx:vir}. We specify the construction because positivity and the topology of the base matter here.

Take the category $\mathcal C$ of finitely presented additive covariant functors $\mathcal D\to\mathrm{Ab}$. It is abelian with pointwise exactness: $\mathcal D$ has cokernels, so $\mathcal D^{\mathrm{op}}$ has weak kernels and Freyd's functor-category criterion applies~\cite[Theorem~3.1]{sx:kv}. Every $F$ has a presentation $H_Y\to H_X\to F\to0$. By Yoneda and left exactness of $\mathrm{Hom}$, the inducing arrow $u:X\to Y$ gives a projective resolution
\[
 0\to H_{\mathrm{coker}\,u}\to H_Y\to H_X\to F\to0.
\]
Every projective in $\mathcal C$ is representable since a summand of a finite sum of representables corresponds to a split idempotent in $\mathcal D$. Define $L(F)=\rho(X)-\rho(Y)+\rho(\mathrm{coker}\,u)$. It is nonnegative by the upper rank inequality. Its value is independent of the resolution by the Euler characteristic identity for finite projective resolutions, since $H_X\mapsto\rho(X)$ is additive on projectives. The horseshoe lemma gives additivity in short exact sequences. Conversely set $\rho(X)=L(H_X)$. The resolution gives the upper inequality and $H_Z\hookrightarrow H_Y$ for an epimorphism $Y\to Z$ gives the lower. These constructions are inverse.

Take $U=H_R$. Every $H_X$ embeds into a finite sum of $U$ by presenting $X$ as a quotient of a finite free module. Every $F$ is a quotient of some $H_X$, so Lemma~\ref{sx:length} applies. Normalization is $L(U)=1$. Each coordinate of every identification and its inverse uses finitely many evaluations, giving an affine homeomorphism with $P(R)$.
\end{proof}

\subsection{Comparing Hilbert modules before taking ranks}
For positive $a\in M_k(A)$ put $E_a=\overline{aA^k}$. For a Hilbert $A$-module $E$, $\mathcal K(E)$ denotes its algebra of compact endomorphisms, with $\theta_{v,w}x=v\langle w,x\rangle$. A submodule is compactly contained in $E$ if a self-adjoint compact endomorphism of $E$ acts as its identity. We call an $A$-linear, inner-product-preserving map an isometric embedding.

\begin{lemma}\label{sx:buffer}
Let $E,F$ be countably generated Hilbert $A$-modules. Suppose every compactly contained submodule of $E$ embeds isometrically in $F$. For each integer $n\ge1$ there is an isometric embedding $E^{\oplus n}\to F^{\oplus(n+1)}$.
\end{lemma}
\begin{proof}
For $z$ in any unital $C^*$-algebra $B$, put $Z=\mathrm{diag}(z,\ldots,z,0)\in M_{n+1}(B)$. Let $S$ be the cyclic shift with ones on the superdiagonal and in the lower left corner. For $s>0$, $Z+sS$ is invertible: solve $(Z+sS)u=v$ using $u_1=s^{-1}v_{n+1}$ followed by $u_{j+1}=s^{-1}(v_j-zu_j)$. If $U_s$ is its unitary polar factor, then
\[
 \|U_s|Z|-Z\|\le\||Z|-|Z+sS|\|+s\longrightarrow0.
\]

Choose a strictly positive $h\in\mathcal K(E)$, omitting the trivial case $E=0$. The submodules $E_j=\overline{(h-1/j)_+E}$ are compactly contained by functional calculus, increase to a dense union, and are countably generated as closed images of a countable family of module generators. Choose embeddings $\phi_j:E_j\to F$ and a countable list of module generators from the union with dense $A$-linear span in $E$. At stage $j$ control the first $j$ listed vectors that already lie in $E_j$. A positive finite-rank contraction $t=\sum_i\theta_{v_i,v_i}$ in $\mathcal K(E_j)$ may be chosen arbitrarily close to the identity on this list, by a contractive approximate identity and approximation of its square root by finite-rank operators. A positive operator of the form $k^*k$ with $k$ finite rank is a finite sum of $\theta_{v,v}$, by factoring the finite Gram matrix.

Set $\alpha=\phi_j$ and $\beta=\phi_{j+1}|_{E_j}$. For $z=\sum_i\theta_{\alpha v_i,\beta v_i}\in\mathcal K(F)$ we have $z\beta x=\alpha tx$ and $|z|\beta x=\beta tx$. Indeed $z^*z$ is the image of $t^2$ under the canonical injective $*$-homomorphism $\mathcal K(E_j)\to\mathcal K(F)$ carrying $\theta_{v,w}$ to $\theta_{\beta v,\beta w}$; this map is available for a Hilbert submodule without an adjointable inclusion~\cite{sx:cei}. Take $B=\mathcal K(F)+\mathbb C I_F$. Let $i:F^{\oplus n}\to F^{\oplus(n+1)}$ append zero. For a unitary $U_j$ as above with $\|U_j|Z|-Z\|<\gamma$, insertion of $t$ on each side gives, for $x\in E_j^{\oplus n}$,
\[
 \|U_j i\beta^{\oplus n}x-i\alpha^{\oplus n}x\|
 \le 2\|(I-t^{\oplus n})x\|+\gamma\|x\|.
\]
Choose $t,\gamma$ so this is at most $2^{-j}$ for each controlled vector in each coordinate. Put $W_1=I$ and $W_{j+1}=W_jU_j$. The maps $W_j i\phi_j^{\oplus n}$ are eventually defined and Cauchy on every listed generator-coordinate vector, hence converge on their algebraic $A$-span. The limit preserves inner products since each map does, and extends by continuity to the claimed isometric embedding.
\end{proof}

\subsection{A closed face of the rank space}
\begin{proof}[Proof of Theorem~\ref{sx:main}]
Let $\mathcal F$ consist of those $r\in P(A)$ satisfying $r(a)=r(f(a))$ for each positive finite matrix $a$ and every continuous $f$ on $\operatorname{sp}(a)\cup\{0\}$ that vanishes at zero and is strictly positive elsewhere. For $g(t)=tf(t)$, both $r(a)-r(g(a))$ and $r(f(a))-r(g(a))$ are nonnegative continuous affine functions on $P(A)$ by the product inequality. The intersection of their zero sets over all such $a,f$ is precisely $\mathcal F$, so $\mathcal F$ is a closed face.

Fix $r\in\mathcal F$. If $b\in M_l(A)_+$ and $c\in\overline{bM_l(A)b}$, then $r(c)\le r(b)$. To see this, act on $M_l(A)$ by left multiplication with $C_0(\operatorname{sp}(b)\setminus\{0\})$ through functional calculus. Its closed essential submodule contains $c$. Cohen--Hewitt factorization gives $c=f(b)y$ with $f$ continuous, vanishing at zero, and $y\in M_l(A)$. Put $k(t)=\sqrt{|f(t)|}+t$. The quotient $f/k$ extends continuously at zero with value zero. Thus $r(c)\le r(k(b))=r(b)$.

Next suppose $T:E_a\hookrightarrow E_b$ is an isometric embedding. Form the rectangular matrix $z$ with columns $T(a^{1/2}e_i)$, where the $e_i$ are the coordinate vectors. These vectors are defined since $E_{a^{1/2}}=E_a$. Then $z^*z=a$. Since each column lies in $E_b$, approximation by columns in $bA^l$ shows $zz^*\in\overline{bM_l(A)b}$. For any rectangular $z$, Hung and Li proved $r(z)=r(|z|)=r(z^*)$~\cite{sx:hl}. The defining power invariance in $\mathcal F$ thus gives $r(z^*z)=r(z)=r(zz^*)$. Consequently $r(a)\le r(b)$.

If $a\precsim b$ in the Cuntz order, the Hilbert-module picture of Coward, Elliott and Ivanescu~\cite{sx:cei} embeds each compactly contained submodule of $E_a$ isometrically in $E_b$. The finitely many columns of $a$ and $b$ generate the respective Hilbert modules after closure. Apply Lemma~\ref{sx:buffer} and the preceding implication to diagonal block matrices. We obtain $nr(a)\le(n+1)r(b)$ for every $n\ge1$. Finiteness allows division by $n$ and gives $r(a)\le r(b)$. Restriction to positive matrices therefore defines a dimension function.

Conversely every dimension function $d$ gives a Sylvester rank $r_d(z)=d([z^*z])$ by the Malcolmson-to-Cuntz comparison in~\cite{sx:hl}. Here $a$, $a^2$ and $f(a)$ generate the same hereditary subalgebra, hence have the same Cuntz class. Thus $r_d\in\mathcal F$, and restriction inverts this assignment in both directions since $r(z^*z)=r(z)$ for $r\in\mathcal F$. Both maps are affine and pointwise continuous, as each evaluation on one side is a fixed evaluation on the other. A nonempty closed face of a compact Choquet simplex is a compact Choquet simplex. Proposition~\ref{sx:rank} now proves the assertion, and the existence theorem of Blackadar and Handelman~\cite{bh:original} supplies a dimension function for every nonzero unital stably finite $A$.
\end{proof}

\medskip
\paragraph{Manuscript history.}
The October~30, 2026 Zenodo preprint~\cite{ap:original} records
the direct-sum argument conditional on the renorming input
\eqref{ap:input}. The present revision states the two-sided separation
in Theorem~\ref{ap:main}.
The Bochner-space and free-lattice
arguments develop~\cite{masw,mafree}. Sources for the mixed-norm and
dimension-function results are given in the corresponding sections.
AI tools assisted with mathematical exploration, literature searches,
and preparation of this manuscript.

\end{document}